\documentclass{amsart}
\usepackage[utf8]{inputenc}
\usepackage{amsmath,amsthm}
\usepackage{comment}
\usepackage{amssymb}
\usepackage{stmaryrd}
\usepackage{hyperref}
\usepackage{color}
\usepackage[left=3cm,right=3cm,top=3cm,bottom=4.5cm]{geometry}
\usepackage{pdfsync}
\usepackage{tikz-cd}
\usepackage{url}
\usepackage{enumitem}
\hypersetup{linktoc=all}
\usepackage[capitalise,noabbrev]{cleveref}
\usepackage{booktabs}

\newcommand{\defterm}[1]{\emph{#1}}

\theoremstyle{plain}
	\newtheorem{thm}{Theorem}
	\numberwithin{thm}{section}
	\newtheorem*{thm*}{Theorem}
	\newtheorem{cor}[thm]{Corollary}
	\newtheorem*{cor*}{Corollary}
	\newtheorem{prop}[thm]{Proposition}
	\newtheorem*{prop*}{Proposition}
	\newtheorem{lem}[thm]{Lemma}
	\newtheorem*{lem*}{Lemma}
	
	\newtheorem*{ex*}{Exercise}
	
	\newtheorem*{claim*}{Claim}
	
	\newtheorem*{question*}{Question}
	\newtheorem{fact}[thm]{Fact}
	\newtheorem*{fact*}{Fact}
\theoremstyle{definition}
	\newtheorem{Def}[thm]{Definition}
	\newtheorem*{Def*}{Definition}
	\newtheorem{obs}[thm]{Observation}
	\newtheorem*{obs*}{Observation}
	\newtheorem{rmk}[thm]{Remark}
	\newtheorem*{rmk*}{Remark}
	
	\newtheorem{soln*}{Solution}
	
	\newtheorem*{note*}{Note}
	\newtheorem{eg}[thm]{Example}
	\newtheorem*{eg*}{Example}	
	
	\newtheorem*{construction*}{Construction}
	
	\newtheorem*{warning*}{Warning}
	
	\newtheorem*{conj*}{Conjecture}
	
    \newtheoremstyle{component}{}{}{}{}{\itshape}{.}{.5em}{\thmnote{#3}#1}
    \theoremstyle{component}
    \newtheorem*{component}{}
	
\newcommand{\nats}{\mathbb{N}}

\newcommand{\op}{\mathrm{op}}
\newcommand{\id}{\mathrm{id}}

\newcommand{\Hom}{\mathrm{Hom}}

\newcommand{\Ob}{\operatorname{Ob}}

\newcommand{\inv}{^{-1}}

\newcommand{\Fun}{\mathrm{Fun}}

\newcommand{\Set}{\mathsf{Set}}

\newcommand{\calC}{\mathcal{C}}

\newcommand{\cube}[1]{\square^{#1}}

\newcommand{\Cubez}{\square_\ast}

\newcommand{\Cube}{\square}
\newcommand{\Cubee}{\square_s}
\newcommand{\Cuber}{\square_r}

\newcommand{\Cubecm}{\square_{c^\wedge}}
\newcommand{\Cubecp}{\square_{c^\vee}}
\newcommand{\Cubecme}{\square_{c^\wedge,s}}
\newcommand{\Cubecpe}{\square_{c^\vee,s}}
\newcommand{\Cubec}{\square_c}
\newcommand{\Cubecr}{\square_{c,r}}
\newcommand{\Cubece}{\square_{c,s}}
\newcommand{\Cubecer}{\square_{c,s,r}}
\newcommand{\Cubeds}{\square_{d,s}}
\newcommand{\Cubedsr}{\square_{d,s,r}}
\newcommand{\Cubedcms}{\square_{d,c^\wedge,s}}
\newcommand{\Cubedcps}{\square_{d,c^\vee,s}}
\newcommand{\Cubedcs}{\square_{d,c,s}}
\newcommand{\Cubedcsr}{\square_{d,c,s,r}}

\newcommand{\cset}{\widehat{\Cube}}

\newcommand{\ccset}{\widehat{\Cubec}}

\newcommand{\cceset}{\widehat{\Cubece}}
\newcommand{\ccerset}{\widehat{\Cubecer}}

\newcommand{\Cubea}{\Cube_{\bullet}}

\newcommand{\cseta}{\widehat{\Cubea}}

\DeclareFontFamily{U}{min}{}
\DeclareFontShape{U}{min}{m}{n}{<-> udmj30}{}
\newcommand\yon{\!\text{\usefont{U}{min}{m}{n}\symbol{'207}}\!}

\tikzset{
  symbol/.style={
    draw=none,
    every to/.append style={
      edge node={node [sloped, allow upside down, auto=false]{$#1$}}}
  }
}

\newcommand{\FinDist}{\mathsf{FinDist}}
\newcommand{\FinLat}{\mathsf{FinLat}}
\newcommand{\Fin}{\mathsf{Fin}}

\newcommand{\fundproj}{\pi}
\newcommand{\fundsym}{\sigma}

\title{Cubical sites as Eilenberg-Zilber categories}
\author{Timothy Campion}
\date{\today}

\begin{document}
\begin{abstract}
    We show that various cube categories (without diagonals, but with symmetries / connections / reversals) are Eilenberg-Zilber categories. This generalizes a result of Isaacson for one particular cubical site. Our method does not involve direct verification of any absolute pushout diagrams. While we are at it, we record some folklore descriptions of cube categories with diagonals and determine exactly which of these are EZ categories.
    
    Beforehand, we develop some general theory of Eilenberg-Zilber categories. We show that a mild generalization of the EZ categories of Berger and Moerdijk are in fact characterized (among a broad class of ``generalized Reedy categories") by the satisfaction of the Eilenberg-Zilber lemma, generalizing a theorem of Bergner and Rezk in the strict Reedy case. We also introduce a mild strengthening of Cisinski's notion of a \emph{cat\'egorie squelettique}, and show that any such category satisfies the Eilenberg-Zilber lemma. It is this tool which allows us to avoid checking absolute pushouts by hand.
\end{abstract}

\maketitle

\tableofcontents


\section{Introduction}
\begin{figure}
    \centering
    \begin{tikzcd}
        \text{Reedy Categories \cite{reedy}} \ar[r,symbol=\subseteq] & 
        \substack{\text{Generalized Reedy Categories} \\ \text{(Def \ref{def:gen})}} \ar[r,symbol=\supseteq] \ar[dr,phantom,"\ulcorner"{description,very near end},"\text{Prop \ref{prop:comp}}" description ] &
        \substack{\text{Generalized Reedy Categories} \\ \text{with right class the monomorphisms}} \\
        \text{Elegant Reedy Categories \cite{bergner-rezk}} \ar[u,symbol=\subseteq] \ar[r,symbol=\subseteq] &
        \substack{\text{Eilenberg-Zilber Categories} \\ \text{(Def \ref{def:ez})}} \ar[u,symbol = \subseteq] \ar[r,symbol=\supseteq]  &
        \text{EZ Reedy Categories \cite{berger-moerdijk}} \ar[u,symbol=\subseteq] \\
        \substack{\text{Reedy Categories satisfying} \\ \text{the Eilenberg-Zilber lemma}} \ar[u,equal,"\text{\cite{bergner-rezk}}"] \ar[r,symbol=\subseteq] &
        \substack{\text{Generalized Reedy Categories satisfying} \\ \text{the Eilenberg-Zilber lemma (Def \ref{def:ez-lemma})}} \ar[u,equal,"\text{Prop \ref{prop:mainprop}}"] \ar[ur,symbol=\supseteq,"\text{\cite{berger-moerdijk}}" below] \\
        \substack{\text{Skeletal} \\ \text{Reedy Categories \cite{cisinski}}} \ar[u,symbol=\subseteq,"\text{\cite{cisinski}} ~"] \ar[r,symbol=\subseteq] &
        \substack{\text{Strongly Skeletal}\\ \text{Generalized Reedy Categories (Def \ref{def:skel})}} \ar[u,symbol=\subseteq,"\text{Thm \ref{thm:ssgrEZ}}~"]
    \end{tikzcd}
    \caption{The types of categories appearing in this note. Inclusions are denoted by ``$\subseteq$", and some inclusions are equalities, denoted with a long ``$=$". References are given for new definitions and implications (mostly in the central column) and citations for definitions and non-obvious inclusions and equalities which are not new. All squares in the diagram are pullbacks; this is by definition except in the case where the reference to Prop \ref{prop:comp} is given. }
    \label{fig:overview}
\end{figure}
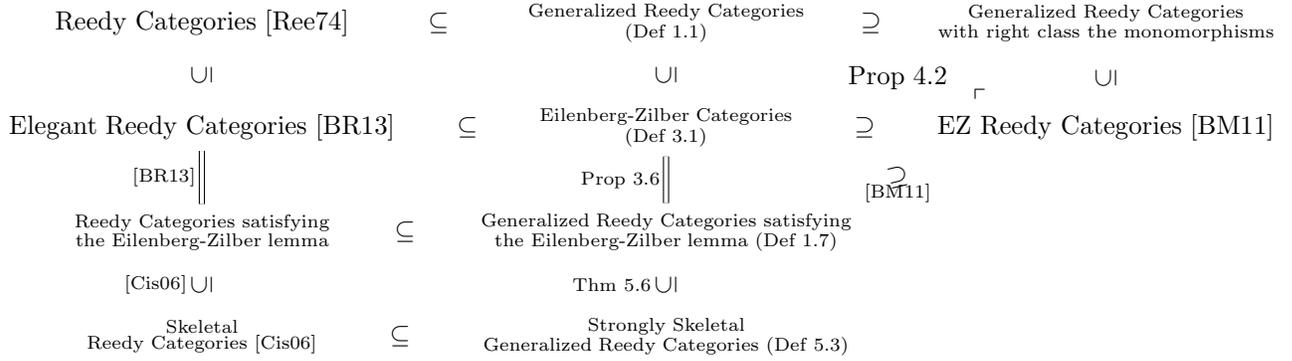

The foundational \emph{Eilenberg-Zilber lemma} \cite{eilenberg-zilber}, which says that every simplex of a simplicial set is -- in a unique way! -- a degeneracy of a nondegenerate simplex, has been generalized in several ways to other indexing categories. For one, the \emph{elegant Reedy categories} of Bergner and Rezk \cite{bergner-rezk} are precisely the Reedy categories whose presheaves satisfy the Eilenberg-Zilber lemma. For another, the \emph{EZ categories} of Berger and Moerdijk \cite{berger-moerdijk} are a type of generalized Reedy category whose presheaves likewise satisfy the Eilenberg-Zilber lemma. In this short note, we show that a mild generalization of EZ categories, which we call \defterm{Eilenberg-Zilber categories} (\cref{def:ez}) are in fact characterized (among a broad class of ``generalized Reedy categories") by the satisfaction of the Eilenberg-Zilber lemma. The notion of \cite{berger-moerdijk} differs only in that it additionally requires the right half of the relevant factorization system to coincide with the class of monomorphisms (\cref{prop:comp}). For the experts, these results should not come as a surprise; rather the surprise may lie in remembering that Berger and Moerdijk did indeed additionally require a monomorphism condition in their work. This work should be seen as a routine verification that the monomorphism assumption is not necessary.
    
We also introduce (\cref{def:skel}) a mild strengthening of Cisinski's notion of a \emph{cat\'egorie squelettique} \cite[Def 8.1.1]{cisinski}, which we call a \defterm{strongly skeletal generalized Reedy category}. We show that any such category satisfies the Eilenberg-Zilber lemma. This generalizes a result of \cite{cisinski} from the Reedy case. We do not know whether every cat\'egorie squelettique satisfies the Eilenberg-Zilber lemma.

As an application, we show in \cref{thm:cubereedy} that several categories of cubes are strongly skeletal, and hence EZ. Cubical sets have received renewed interest in higher category theory and homotopy type theory. The sites we treat are recalled in \cref{sec:cubsite}. It was previously shown in \cite{isaacson} that the site of cubes with one connection and symmetries is EZ; here we extend to allow two connections and/or reversals as well. Our use of the strong skeletality property simplifies the proof even in the case of Isaacson's site: unlike Isaacson, we do not need to manually construct any absolute pushouts. Unfortunately, cubical sites with \emph{diagonals}, which are generally required for applications to homotopy type theory, are not generally EZ, as we show in \cref{sec:diag}.\footnote{We will apply our results in \cite{campion-subdiv} in a study of cubical subdivision, with more homotopy-theoretic applications in mind.} 

The outline of the note is as follows. After some basic definitions in \cref{sec:genreedy}, we recall some background on absolute colimits in \cref{sec:absolute}, which we use to define Eilenberg-Zilber categories in \cref{sec:ez}. We compare this definition to the one in \cite{berger-moerdijk} in \cref{sec:bm}. We discuss our modification of Cisinski's skeletality condition in \cref{sec:skel} and its relationship to the Eilenberg-Zilber property. In the final two sections, we give an application to cubical sets. The necessary background is recalled in \cref{sec:cubsite}, and the main theorem, saying that several of these sites are Eilenberg-Zilber appears in \cref{sec:cubsitegenreedy}. In \cref{sec:diag}, we record what happens in the cases where diagonals are allowed. In all such cases, there is either an obstruction in the form of non-idempotent completeness, or else an alternate combinatorial description of the category allowing a direct deduction of the EZ property. In the non-idempotent complete cases, we determine whether the idempotent completion is EZ.

\subsection{Notation} For a category $\calC$ we denote $\widehat \calC = \Fun(\calC^\op,\Set)$ the category of presheaves on $\calC$. We denote by $\yon_\calC : \calC \to \widehat \calC$ the Yoneda embedding.

\subsection{Acknowledgements} I would like to thank Denis-Charles Cisinski, Sina Hazratpour, Chris Kapulkin, Keith Kearnes, Emily Riehl, Maru Sarazola, and Jonathan Weinberger for helpful conversations. I am grateful for the support of NSF grant DMS-1547292 and the hospitality of MSRI's Spring 2020 \emph{Higher Categories and Categorification} program during the writing of this note, as well as the ARO under MURI Grant W911NF-20-1-0082.



\section{Generalized Reedy Categories}\label{sec:genreedy}
We begin with the following ``minimal" definition, including all but one axiom of \cite{berger-moerdijk}:

\begin{Def}\label{def:gen}
A \defterm{generalized Reedy category} is an essentially small category $A$ equipped with 
\begin{itemize}
    \item a ``degree" function $\deg: \Ob A \to \nats$, which respects isomorphisms, and
    \item an orthogonal factorization system $(A_-,A_+)$,
\end{itemize}  
such that 
\begin{itemize}
    \item every non-isomorphism of $A_-$ strictly lowers degree, and
    \item every non-isomorphism of $A_+$ strictly raises degree.
\end{itemize}
(That is, if $(f: a \to b) \in A_-$, then $\deg(a) > \deg(b)$ and iff $(f :a \to b) \in A_+$, then $\deg(a) < \deg(b)$.)
\end{Def}

Here, we recall that 
\begin{Def}
An \defterm{orthogonal factorization system} on a category $A$ consists of two class of morphisms $A_-,A_+$ in $A$ such that
\begin{enumerate}
    \item[(i)] Every morphism admits a factorization as a morphism in $A_-$ followed by a morphism in $A_+$, and
    \item[(ii)] The morphisms of $A_-$ have the unique left lefting property against morphisms of $A_+$.
\end{enumerate}
Recall also that in the presence of (i), condition (ii) is equivalent \cite[Appendix C]{joyal} to
\begin{enumerate}
    \item[(ii')] $A_-$ and $A_+$ each contain the isomorphisms and are closed under composition, and the factorizations of (i) are unique up to unique isomorphism.
\end{enumerate}
\end{Def}

\begin{rmk}
There are several notions of ``generalized Reedy category" in the literature (cf. \cite{berger-moerdijk}, \cite{cisinski}, \cite{shulman}). We make no claims about the usefulness of this definition in general, but merely introduce it for the purposes of comparison in this note.
\end{rmk}

We check that, as in the Reedy case \cite[Lemma 2.9]{riehl-verity}, a generalized Reedy structure in fact depends only on the underlying category and degree function.

\begin{lem}[cf. \cite{riehl-verity}]
Let $A$ be an essentially small category equipped with a ``degree" function $\deg: \Ob A \to \nats$. Then there is at most one generalized Reedy structure on $A$ compatible with this degree function.
\end{lem}
\begin{proof}
Suppose that $A$ admits a generalized Reedy structure using the given degree function. Then $A_-$ may be characterized as the collection of extremally degree-lowering morphisms. That is, if $a \xrightarrow f b$ is a morphism in $A$, then $f \in A_-$ if and only if (1) $\deg b \leq \deg a$ and (2) for any factorization $a \xrightarrow f b = a \xrightarrow g c \xrightarrow h b$, with $\deg c \leq \deg b$, we have that $h$ is an isomorphism. Of course, $A_+$ admits a dual characterization.

For one direction, suppose that $f$ satisfies (1) and (2), and factor $a \xrightarrow f b = a \xrightarrow g c \xrightarrow h b$ where $g \in A_-$ and $h \in A_+$. Then $h$ must weakly raise degree, so that $\deg c \leq \deg b$, and so by (2) we have that $h$ is an isomorphism. Since $(A_-,A_+)$ is an orthogonal factorization system, it is invariant under isomorphism, so because $g \in A_-$ we have $f \in A_-$.

Conversely, suppose that $a \xrightarrow f b \in A_-$. If $f$ is an isomorphism, then $\deg a = \deg b$ because $\deg$ respects isomorphism; if $f$ is not an isomorphism, then $\deg b < \deg a$ by hypothesis. Thus (1) holds. To see that (2) holds, suppose that $a \xrightarrow f b = a \xrightarrow g c \xrightarrow h b$ where $\deg c \leq \deg b$. First consider the case where $a \xrightarrow g c \in A_-$ and $c \xrightarrow h b \in A_+$. In this case, by uniqueness of $(A_-,A_+)$ factorizations we have that $h$ is an isomorphism, as desired. Now consider the case where $a \xrightarrow g c \in A$ is arbitrary and $c \xrightarrow h b \in A_+$. Then we may factor $a \xrightarrow g c = a \xrightarrow {g_1} a_1 \xrightarrow {g_2} c$ where $g_1 \in A_-, g_2 \in A_+$. Then since $hg_2 \in A_+$, we have by the previous case that $hg_2$ is an isomorphism. It follows that $h,g_2$ preserve degree, and therefore $h$ is an isomorphism, as desired. Finally, consider the case where $a \xrightarrow g c \in A$ and $c \xrightarrow h b \in A$ are both arbitrary. Factor $c \xrightarrow h b = c \xrightarrow {h_1} b_1 \xrightarrow {h_2} b$ where $h_1 \in A_-, h_2 \in A_+$. Then because $h_2 \in A_+$, we have by the previous case that $h_2$ is an isomorphism. Therefore $h_1 \in A_+$ is weakly degree-raising and so $h_1$ is an isomorphism. Thus $h = h_2 h_1$ is also an isomorphism as desired.
\end{proof}

Just like Reedy categories \cite{bergner-rezk}, generalized Reedy categories are idempotent-complete:

\begin{lem}[cf. \cite{bergner-rezk}]\label{lem:idems}
Let $A$ be a generalized Reedy category. Then idempotents split in $A$.
\end{lem}
\begin{proof}
Let $e$ be an idempotent, and factor $e =gf$ with $f \in A_-$, $g \in A_+$. Factor $fg = g'f'$ with $f' \in A_-$, $g' \in A_+$. Then from the equation $e = e^2$, we have $gf = gfgf = gg'f'f$. By uniqueness of factorizations, we have up to isomorphism that $g = gg'$, $f = f'f$. Therefore $f',g'$ are degree-preserving and hence isomorphisms since they lie in $A_-$ and $A_+$ respectively. Therefore $fg = g'f'$ is an isomorphism. At the same time, we have $(fg)^3 = f(gf)^2 g = f(gf)g = (fg)^2$. But any isomorphism satisfying this equation is the identity, so $fg = 1$. Thus the factorization $e = gf$ is a splitting of $e$.
\end{proof}

Generalized Reedy categories are stable under slicing, and more generally under passage to categories of elements :

\begin{lem}\label{lem:grslice}
Let $(A,A_-,A_+)$ be a generalized Reedy category, and let $X \in \widehat A$. Then the category of elements $A \downarrow X$ is a generalized Reedy category with the degree function defined as in $A$.
\end{lem}
\begin{proof}
More generally, factorization systems are stable under passage to categories of elements.
\end{proof}

The notion of Eilenberg-Zilber decomposition makes sense in this generality; the Eilenberg-Zilber lemma may be stated, but will not hold in all generalized Reedy categories:

\begin{Def}\label{def:ez-lemma}
Let $A$ be a generalized Reedy category. For $X \in \widehat A$ and $a \in A$, we write $X_a := X(a) = \Hom(a,X)$ for the set of \defterm{$a$-cells} of $X$. An $a$-cell $a \to X$ is \defterm{nondegenerate} if it does not factor through any non-isomorphism $a \to b$ in $A_-$. An \emph{Eilenberg-Zilber decomposition} of an $a$-cell $a \to X$ is a factorization $a \to b \to X$ where $a \to b \in A_-$ and $b \to X$ is nondegenerate. We say that an Elienberg-Zilber decomposition is \emph{essentially unique} if for any two such decompositions $a \to b \to X$, $a \to b' \to X$, there is a unique isomorphism $b \to b'$ such that $a \to b \to b' = a \to b'$ and $b \to b' \to X = b \to X$. We say that $A$ \defterm{satisfies the Eilenberg-Zilber lemma} if for every $a$-cell $a \to X$ for $a \in A, X \in \widehat A$, there is an essentially unique Eilenberg-Zilber decomposition.
\end{Def}

As usual though, the existence part of the Eilenberg-Zilber lemma does hold in any generalized Reedy category:

\begin{lem}\label{lem:ez-exist}
Let $A$ be a generalized Reedy category. Then for every presheaf $X \in \widehat A$ and every cell $x \in X_a$, there exists a factorization of $x = x'f$ through a representable of minimal degree. Moreover, any such factorization is an Eilenberg-Zilber decomposition of $x$ (which is not in general essentially unique).
\end{lem}
\begin{proof}
The first statement follows from the fact that the set of degrees is well-ordered. For the second statement, we first show that if $x = x'f$ is a factorization of $x$ through a representable of minimal degree, then $f \in A_-$. For otherwise we could factor $f = hg$ with $g \in A_-, h \in A_+$. By degree minimality, we have that $h$ is degree-preserving and hence an isomorphism. Because $g \in A_-$, if follows that $f \in A_-$. Now suppose for contradiction that $x = x'f$ is a factorization through a representable of minimal degree which is not an Eilenberg-Zilber decomposition. Then $x'$ is degenerate, so that $x' = x''k$ where $k \in A_-$ is not an isomorphism, and hence strictly degree-lowering. But then $x = x''kf$ is a factorization of $x$ through a representable of lower degree, a contradiction. Thus every factorization of $f$ through a representable of minimal degree is an Eilenberg-Zilber decomposition.
\end{proof}


\section{Absolute Colimits}\label{sec:absolute}

Recall that a colimit diagram in a category $\mathcal C$ is called an \defterm{absolute colimit} if the colimit is preserved by any functor $\mathcal C \to \mathcal D$ to any category $\mathcal D$, or equivalently if it is preserved by the Yoneda embedding $\yon_\calC: \calC \to \widehat \calC$.

As seen in \cite{berger-moerdijk} and \cite{bergner-rezk}, the Eilenberg-Zilber lemma is closely related to absolute pushouts, and absolute (=split) epimorphisms. We review some of their properties here.

\begin{lem}\label{lem:split-iso}
Let $a \xrightarrow f b$ be a split epimorphism in a category $\calC$. If there exists $b \xrightarrow g c$ such that $gf$ is a monomorphism, then $f$ is an isomorphism.
\end{lem}
\begin{proof}
If $gf$ is monomorphism, then $f$ is a monomorphism, and any monic split epic is an isomorphism.
\end{proof}

\begin{lem}\label{lem:split-abs}
Let $a \xrightarrow f b$ be a morphism of a category $\calC$. Then the following are equivalent:
\begin{enumerate}
    \item The pushout of $f$ along itself exists and is absolute.
    \item $f$ is an absolute epimorphism.
    \item The image of $f$ under the Yoneda embedding is an epimorphism.
    \item $f$ is a split epimorphism.
\end{enumerate}
\end{lem}
\begin{proof}
Since split epimorphisms are epimorphisms and are preserved by functors, they are absolute epimorphisms; thus $(4) \Rightarrow (2)$. It is clear that $(2) \Rightarrow (3)$. For $(3) \Rightarrow (4)$, note that if $\yon_\calC(f)$ is an epimorphism, then $\Hom(b,a) \xrightarrow{\Hom(b,f)} \Hom(b,b)$ is surjective. In particular, $\id_b$ is in the image of this map, i.e. there exists $g \in \Hom(b,a)$ such that $fg = \id_b$, i.e. $f$ is a split epimorphism.

$(2) \Rightarrow (1)$ is clear because if $f$ is an absolute epimorphism, then $1 f = 1 f$ is an absolute pushout square. Conversely, suppose that that $(1)$ holds. Because the pushout is preserved by the Yoneda embedding, we have that $\Hom(b \cup_a b, b \cup_a b) = \Hom(b\cup_a b, b) \cup_{\Hom(b\cup_a b, a)} \Hom(b\cup_a b, b)$. In particular, $\id_{b \cup_a b}$ must factor through one of the copies of $b$, i.e. one of the maps $b \to b \cup_a b$ is a split epimorphism (and by symmetry so is the other). Therefore, $\Hom(a,b) \to \Hom(a,b\cup_a b) = \Hom(a,b) \cup_{\Hom(a,a)} \Hom(a,b)$ is surjective. But clearly for any map $\phi$ of sets, if the cobase change of $\phi$ along itself is surjective, then $\phi$ is already surjective. Thus $\Hom(a,a) \to \Hom(a,b)$ is surjective, i.e. $f$ is a split epimorphism, so that $(1) \Rightarrow (4)$.
\end{proof}

\begin{cor}\label{cor:abs-split}
Let $a\xrightarrow f b$ be a split epimorphism in a category $\calC$, and suppose that the pushout of $f$ along $a \xrightarrow g c$ exists and is absolute. Then the cobase-change $c \xrightarrow {f'} b \cup_a c$ of $f$ along $g$ is also a split epimorphism.
\end{cor}
\begin{proof}
If $a \xrightarrow f b$ is a split epimorphism, then $\yon_\calC(f)$ is an epimorphism. Since the absolute pushout is preserved by $\yon_\calC$ and epimorphisms are stable under cobase-change, we have that $\yon_\calC(f')$ is an epimorphism, and thus $f'$ is a split epimorphism by Lemma \ref{lem:split-abs}.
\end{proof}


\section{Eilenberg-Zilber Categories}\label{sec:ez}

Inspired by \cite{berger-moerdijk}, we make the following definition. A precise comparison to the EZ categories of \cite{berger-moerdijk} appears in \cref{sec:bm}.

\begin{Def}\label{def:ez}
An \defterm{Eilenberg-Zilber category} is a generalized Reedy category such that every pair of morphisms of $A_-$ with a common domain has an absolute pushout in $A$.
\end{Def}

\begin{lem}\label{lem:ez-split}
Let $A$ be an Eilenberg-Zilber category. Then the morphisms of $A_-$ are precisely the split epimorphisms.
\end{lem}
\begin{proof}
That $A_-$ is contained in the split epimorphisms follows from \ref{lem:split-abs}. For the converse, first note that $A_-,A_+$ have no non-identity idempotents, since an idempotent must be degree-preserving and hence an isomorphism if it is in $A_-$ or $A_+$. Next, if $f \in A_+$ has a section $s \in A$, then factor $s = hg$ with $g \in A_-$, $h \in A_+$. Then $1 = (fh)g$ is an $(A_-,A_+)$ factorization; by uniqueness we may assume that $g =1$ and $fh = 1$. If $f$ is not an isomorphism, then $hf$ is a non-identity idempotent in $A_+$, a contradiction; thus $f$ is an isomorphism. Finally, if $f \in A$ has a section $s$, then factor $f = hg$ with $g \in A_-, h \in A_+$. Then $gs$ is a section of $h \in A_+$, so by the previous step we have that $h$ is an isomorphism; since $g \in A_-$ it follows that $f \in A_-$ as desired.
\end{proof}

\begin{rmk}
It follows from \cref{lem:ez-split} that a category can admit at most one Eilenberg-Zilber structure, up to modifications of the degree function which do not affect the factorization system.
\end{rmk}

\begin{cor}\label{cor:ez-po}
Let $A$ be an Eilenberg-Zilber category, and let $g'f = f'g$ be an absolute pushout of $f,g$ where $f,g \in A_-$. Then $f',g' \in A_-$ as well.
\end{cor}
\begin{proof}
This follows from Lemma \ref{lem:ez-split} and Corollary \ref{cor:abs-split}.
\end{proof}

\begin{rmk}
It follows from Corollary \ref{cor:ez-po} that in Definition \ref{def:ez}, we could equivalently ask that $A_-$ has all pushouts, which are preserved by the inclusion into $A$, and are absolute as pushouts in $A$.
\end{rmk}

The point of Definition \ref{def:ez} is to ensure that the Eilenberg-Zilber lemma is satisfied. Not only is this the case (cf. \cite{berger-moerdijk}), but the converse also holds, as shown in the Reedy case by \cite{bergner-rezk}; the following proof is the same:

\begin{prop}[cf. \cite{berger-moerdijk},\cite{bergner-rezk}]\label{prop:mainprop}
Let $A$ be a generalized Reedy category (\cref{def:gen}). Then $A$ satisfies the Eilenberg-Zilber lemma (\cref{def:ez-lemma}) if and only if $A$ is an Eilenberg-Zilber category (\cref{def:ez}).
\end{prop}
\begin{proof}
For one direction, suppose that $A$ is an Eilenberg-Zilber category. We have seen in Lemma \ref{lem:ez-exist} that Eilenberg-Zilber decompositions always exist, so we need only show that they are essentially unique. Suppose that $xf = yg$ are two Eilenberg-Zilber decompositions of the same $a$-cell. Because $A$ is an Eilenberg-Zilber category, we have an absolute pushout square $g'f = f'g$, and we have a factorization $z$ such that $x = zg'$, $y = zf'$. Moreover, by Corollary \ref{cor:ez-po}, we have $g',f' \in A_-$. But then because $xf$ is an Eilenberg-Zilber decomposition factoring through the decomposition $z(g'f)$, we must have that $g'$ is an isomorphism. Similarly, $f'$ is an isomorphism. Then the isomorphism $(f')\inv g'$ exhibits the decompositions $xf$ and $yg$ as essentially equivalent. Uniqueness of this isomorphism follows from the fact that $f,g$ are epimorphisms.

For the converse, suppose that $A$ satisfies the Eilenberg-Zilber lemma. Let $f,g \in A_-$ have common domain, and let $P$ be their pushout in $\widehat A$, with pushout square $g'f = f'g$. We take Eilenberg-Zilber decompositions of $f',g'$; by composition with $g,f$ these yield Eilenberg-Zilber decompositions of $f'g = g'f$, and by uniqueness we may assume without loss of generality that these decompositions factor through the same nondegenerate cell $c \to P$. By pushout we obtain a map $P \to c$, and it is straightforward to verify that $P \to c \to P$ is the identity. Thus $P$ is a retract of $c$, and by Lemma \ref{lem:idems} we have $P \in A$. Because the Yoneda embedding is fully faithful, the pushout square lies already in $A$, and is a pushout preserved by the Yoneda embedding. Thus $A$ is an Eilenberg-Zilber category.
\end{proof}


\section{Comparison to Berger-Moerdijk}\label{sec:bm}

Berger and Moerdijk actually give a slightly different definition than Definition \ref{def:ez}, which we will show is strictly more restrictive:
\begin{Def}[\cite{berger-moerdijk}]
A \defterm{Berger-Moerdijk EZ category} is a small category equipped with a degree function $\deg: \Ob A \to \nats$ such that
\begin{enumerate}
    \item A monomorphism preserves (resp. raises) degree iff it is invertible (resp. noninvertible);
    \item Every morphism admits a factorization as a split epimorphism followed by a monomorphism;
    \item Every pair of split epimorphisms with a common domain has an absolute pushout.
\end{enumerate}
\end{Def}

\begin{prop}\label{prop:comp}
Every Berger-Moerdijk EZ category is an Eilenberg-Zilber category with the same degree function. Conversely, an Eilenberg-Zilber category $A$ is a Berger-Moerdijk EZ category if and only if every morphism of $A_+$ is a monomorphism.
\end{prop}
\begin{proof}
First, let $A$ be a Berger-Moerdijk EZ category and let us show that $A$ is an Eilenberg-Zilber category, with $A_-$ being the split epimorphisms and $A_+$ being the monomorphisms.

First we show that a split epimorphism preserves (resp. strictly lowers) degree iff it is invertible (resp. noninvertible). If $a \xrightarrow f b$ is a split epimorphism, then let $b \xrightarrow s a$ be a section. Since $s$ is a split monomorphism, we have $s \in A_+$. We have that $f$ preserves (resp. strictly lowers) degree iff $s$ preserves (resp. strictly raises) degree iff $s$ is invertible (resp. noninvertible) iff $f$ is invertible (resp. noninvertible).

Now we show that the degree function $\deg: \Ob A \to \nats$ respects isomorphisms. Let $a \xrightarrow f b$ be an isomorphism, and suppose that $\deg a \neq \deg b$; without loss of generality we have $\deg a < \deg b$. Factor $a \xrightarrow f b = a \xrightarrow g c \xrightarrow h b$ where $g$ is a split epimorphism and $h \in A_+$. By Lemma \ref{lem:split-iso}, $g$ is an isomorphism. Since $g \in A_-$, by the previous part we have $\deg c = \deg a$. Thus $h \in A_-$ strictly raises degree, and so is not an isomorphism. Since $f = hg$ and $g$ is an isomorphism, this implies that $f$ is not an isomorphism, a contradiction.

Now we show that $(A_-,A_+)$ factorizations are unique up to isomorphism. Suppose that $f = hg = h'g'$ where $g,g' \in A_-$ and $h,h' \in A_+$. There is an absolute pushout square $kg = k'g'$, and by Corollary \ref{cor:abs-split}, we have that $k,k'$ are split epimorphisms; it will suffice to show that $k,k'$ are isomorphisms. There is also an induced map $h''$ such that $h'' k = h, h'' k' = h'$.  Factoring $h'' = ml$ with $l \in A_-$, $m \in A_+$, we obtain a factorization $h = mlk$. Because $h$ is a monomorphism, we may invoke Lemma \ref{lem:split-iso} to see that $k$ is an isomorphism as desired.

In fact, $(A_-,A_+)$ factorizations are unique up to unique isomorphism; this follows either from $A_-$ consisting of epimorphisms or $A_+$ consisting of monomorphisms. Since $A_-$ and $A_+$ are also closed under composition, it follows that $(A_-,A_+)$ is an orthogonal factorization system, and $A$ is a generalized Reedy category. It is an EZ category by virtue of (3).

Now let $A$ be an Eilenberg-Zilber category such that every morphism of $A_+$ is a monomorphism. Suppose that $a \xrightarrow f b $ is a monomorphism, and factor $a \xrightarrow f b = a \xrightarrow g c \xrightarrow h b$ where $g \in A_-$, $h \in A_+$. Then by Lemma \ref{lem:split-iso}, $g$ is an isomorphism. It follows that $f \in A_+$. Thus $A_+$ consists exactly of the monomorphisms. By Lemma \ref{lem:ez-split}, $A_-$
consists of exactly the split epimorphisms.
It is now straightforward to see that $A$ is a Berger-Moerdijk EZ category.
\end{proof}

\begin{eg}
Let $A$ be a direct category where not every morphism is a monomorphism. Then $A$ is an Eilenberg-Zilber category, but not a Berger-Moerdijk EZ category. We do not know of other examples.
\end{eg}


\section{Strongly Skeletal Generalized Reedy Categories}\label{sec:skel}
In this section, we introduce \emph{strongly skeletal} generalized Reedy categories (\cref{def:skel}) and prove that they are Eilenberg-Zilber categories (\cref{thm:ssgrEZ}). These are inspired by the \emph{skeletal} generalized Reedy categories of \cite{cisinski}, whose definition we first recall.
\begin{Def}[\cite{cisinski}]
We recall some terminology from \cite{cisinski}:
\begin{itemize}
    \item A \defterm{skeletal generalized Reedy category} (\cite[Def 8.1.1]{cisinski} ``cat\'egorie squelettique") is a generalized Reedy category $(A,A_-,A_+)$ such that (1) each morphism of $A_-$ is a split epimorphism in $A$, and (2) no two distinct morphisms of $A_-$ have the same set of sections in $A$.
    \item A \defterm{skeletal Reedy category}\footnote{This agrees with \cite{cisinski}'s \emph{cat\'egorie squelettique normale} by \cite[Prop 8.1.37]{cisinski}.}  is a skeletal generalized Reedy category which is a Reedy category.
    \item A \defterm{regular skeletal Reedy category} (\cite[Def 8.2.3]{cisinski}) is a skeletal Reedy category where every morphism of $A_+$ is a monomorphism.
\end{itemize}
\end{Def}

Although \cite{cisinski} gives this definition in the generalized setting (i.e. allowing for nontrivial automorphisms), it seems that not much is known about skeletal generalized Reedy categories except in the non-generalized case. Indeed, it seems that in order to use the notion nontrivially, the definition must be strengthened, replacing sets of sections with sets of \emph{pseudo-sections}, as we now define.

\begin{Def}\label{def:pseudo}
If $p: a \to b$ is a morphism in a category, say that $i: c \to a$ is a \defterm{pseudo-section} of $p$ if $pi: c \to b$ is an isomorphism. We say that two maps $p: a \to b$, $p': a \to b'$ are \defterm{pseudo-equal} if there is an isomorphism $\phi: b \to b'$ such that $\phi p = p'$.
\end{Def}

\begin{Def}\label{def:skel}
A \defterm{strongly skeletal generalized Reedy category} is a generalized Reedy category $(A,A_-,A_+)$ such that (1) every morphism of $A_-$ is a split epimorphism in $A$ and (2) for every two morphisms $p: a \to b$, $p': a \to b'$ of $A_-$ with the same domain, if $p,p'$ have identical sets of pseudo-sections in $A$, then $p,p'$ are pseudo-equal.
\end{Def}

\begin{rmk}
It is straightforward to show that any strongly skeletal generalized Reedy category $(A,A_-,A_+)$ is a skeletal generalized Reedy category (in fact, if $p,p': a \to b$ are morphisms in any category with common domain and codomain, and if $p,p'$ have the same set of sections, then they also have the same set of pseudo-sections). If $(A,A_-,A_+)$ is in fact Reedy, then the converse holds because $A$ has no non-identity isomorphisms. Thus any skeletal Reedy category is strongly skeletal. The author does not know an example of a skeletal generalized Reedy category which is not strongly skeletal, but suspects that they exist.
\end{rmk}

\begin{lem}[{cf. \cite[8.1.6]{cisinski}}]\label{lem:ssgrslice}
Let $(A,A_-,A_+)$ be a strongly skeletal generalized Reedy category, and let $X \in \widehat A$. Then the slice category $A \downarrow X$ is a strongly skeletal generalized Reedy category with the degree function defined as in $A$.
\end{lem}
\begin{proof}
We know that $A \downarrow X$ is a generalized Reedy category from Lemma \ref{lem:grslice}. Observe that the set of pseudo-sections of a morphism in $A\downarrow X$ is canonically identified with the set of pseudo-sections of the underlying morphism in $A$. Moreover, observe that if $\pi,\pi'$ are epimorphisms in a category which are pseudo-equal, then they are also pseudo-equal in any slice category.
\end{proof}

\begin{thm}[{cf. \cite[8.1.24]{cisinski}}]\label{thm:ssgrEZ}
Let $(A,A_-,A_+)$ be a strongly skeletal generalized Reedy category. Then $A$ is an Eilenberg-Zilber category.
\end{thm}
\begin{proof}
Let $\sigma: a \to X$ be a cell, and let $a \to b \to X$, $a \to b' \to X$ be two Eilenberg-Zilber decompositions; we would like to show they are equivalent (recall from Lemma \ref{lem:ez-exist} that some Eilenberg-Zilber decomposition always exists). By Lemma \ref{lem:ssgrslice}, we may pass to $A\downarrow X$ and thus assume that $X$ is the terminal object. With this assumption, $b$ and $b'$ become nondegenerate objects, meaning that every map $b \to c$ or $b' \to c$ in $A_-$ is an isomorphism. Because $a \to b$, $a \to b'$ are in $A_-$, they have sections. By composition, we obtain maps $b \to b'$ and $b' \to b$ which must be in $A_+$ by nondegeneracy of $b,b'$. Since both maps are degree-non-decreasing, it follows that $b,b'$ have the same degree. It then follows that the maps $b \to b'$, $b' \to b$ are isomorphisms. So we may assume without loss of generality that $b = b'$. 

Thus it suffices to show that if $A$ is a strongly skeletal generalized Reedy category and $b \in A$ is a nondegenerate object, then any two maps $p,p': a \rightrightarrows b$ in $A_-$ are related by an automorphism under $a$, i.e. that they are pseudo-equal (because $p,p'$ are epimorphisms, the automorphism is automatically unique). By strong skeletality, it suffices to show that $p,p'$ have the same sets of pseudo-sections. In fact, we claim that $b'' \to a$ is a pseudo-section of $b$ iff it is a pseudo-section of $a$, iff there exists an isomorphism $b'' \cong b$. In one direction, for any pseudo-section $b'' \to a$ of $p$ or $p'$, we have $b\cong b''$. Conversely, suppose that $b'' \cong b$. Then $b''$ is nondegenerate and has the same degree as $b$, so any map $b'' \to b$ is an isomorphism. Consequently, every map $b'' \to a$ is a pseudo-section of any map $a \to b$.
\end{proof}




\section{Cubical sites}\label{sec:cubsite}
We recall from \cite{grandis} and \cite{buchholtz-morehouse} several categories of ``cubes", and the associated presheaf categories of ``cubical sets". Terminology and notation varies; ours most closely resembles that of Isaacson \cite{isaacson}. An exhaustive taxonomy of cubical sites, using the framework of monoidal or substructural logic of \cite{mauri}, is given in \cite{buchholtz-morehouse}. We are most interested in the sites considered in \cite{grandis}, but we consider more generally cubical sites with \emph{weakening} and without \emph{contraction} in the terminology of \cite{buchholtz-morehouse}. 

Let $\Cubez$ be the full subcategory of $\Set$ on those sets of the form $\{0,1\}^n$ where $n$ is finite. We freely identify $\{0,1\}^n$ with the powerset lattice $P(\{1,\dots,n\})$ via the bijection $(\epsilon_1,\dots, \epsilon_n) \mapsto \{i \mid \epsilon_i = 1\}$, and denote it $\cube n$. Note that $\Cubez$ is (cartesian) symmetric monoidal, with $\cube m \otimes \cube n = \cube {m+n}$. All of our categories of cubes will be non-full, noncartesian, generally nonsymmetric, wide monoidal subcategories of $\Cubez$.

\textbf{Plain cubical sets.} Let $\Cube$ be the wide monoidal subcategory of $\Cubez$ generated by the two \defterm{fundamental face maps} $\iota^0,\iota^1: \cube 0 {}^\to_\to \cube 1$ and the \defterm{fundamental projection map} $\fundproj:\cube 1 \to \cube 0$. So a map $\cube m \to \cube n$ is an $n$-fold tensor product of maps $\cube {m_i} \to \cube 1$ where $m_1 + \dots + m_n = m$, and a map $\cube {m_i} \to \cube 1$ is either constant at 0 or 1, or else is a projection onto some factor. Equivalently (\cite[8.4.5]{cisinski}, \cite{grandis}), $\Cube$ is the free monoidal category on the two-object category generated by the fundamental face and projection maps. The category of \defterm{plain cubical sets} is the presheaf category $\cset$.

\textbf{Cubical sets with connections.} Let $\Cubec$ be the wide monoidal subcategory of $\Cubez$ generated by $\Cube$ along with the \defterm{fundamental connection maps} $\gamma_\vee, \gamma_\wedge: \cube 2 {}^\to_\to \cube 1$. Here $\gamma_\vee$ is union and $\gamma_\wedge$ is intersection in the lattice $\cube 1 = P(\{1\})$. Again a map $\cube m \to \cube n$ is an $n$-fold tensor product of maps $\cube {m_i} \to \cube 1$. Note that $\gamma_\vee$ (resp. $\gamma_\wedge$) is associative with unit $\iota^0$ (resp. $\iota^1$) and absorbing element $\iota^1$ (resp. $\iota^0$); such a structure is called a \defterm{cubical monoid}. Equivalently \cite{grandis}, $\Cubec$ is the free monoidal category generated by a cubical monoid. The category of \defterm{cubical sets with connections} is the presheaf category $\ccset$.

There are also subcategories $\Cubecm,\Cubecp$ which have only one type of connection (and are of course equivalent, but not under $\Cube$). They are studied in \cite{maltsiniotis-cube},\cite{cisinski-elegant},\cite{isaacson},\cite{buchholtz-morehouse},\cite{kapulkin-lindsey-wong} and elsewhere.

\textbf{Cubical sets with connections and symmetries.} Let $\Cubece$ be the wide monoidal subcategory of $\Cubez$ generated by $\Cubec$ along with the \defterm{fundamental symmetry map} (also called an \defterm{exchange} in \cite{buchholtz-morehouse} or \defterm{extension} map in \cite{grandis}) $\fundsym: \cube 2 \to \cube 2$, $(x,y) \mapsto (y,x)$. Equivalently \cite{grandis}, $\Cubece$ is the free monoidal category on a \defterm{symmetric cubical monoid} (i.e. a cubical monoid $M$ with an involution $\fundsym$ on $M^{\otimes 2}$ such that $\gamma_\pm \fundsym = \gamma_\pm$). Alternatively, it follows from the work of \cite{grandis} that $\Cubece$ is the symmetric monoidal subcategory of $\Cubez$ generated by $\Cubec$, or the free symmetric monoidal category on a symmetric cubical monoid (where now $\fundsym$ is taken to be the symmetry isomorphism of the symmetric monoidal category). The category of \defterm{cubical sets with connections and symmetries} is the presheaf category $\cceset$.

There are also subcategories $\Cubee, \Cubecme,\Cubecpe$ which have symmetries but not connections, or only one type of connection. The sites $\Cubecme,\Cubecpe$ (which are of course equivalent) are studied in \cite{isaacson}.

\textbf{Cubical sets with connections, symmetries, and reversals.} Let $\Cubecer$ be the wide monoidal subcategory of $\Cubez$ generated by $\Cubece$ along with the \defterm{fundamental reversal map} $\rho: \cube 1 \to \cube 1$, $0 \mapsto 1$, $1 \mapsto 0$. Equivalently \cite{mauri}, $\Cubecer$ is the free monoidal category on a \defterm{symmetric involutive cubical monoid} (i.e. a symmetric cubical monoid $M$ equipped with an involution $\rho$ on $M$ such that $\rho \iota^0 = \iota^1$, $\fundproj \rho = \fundproj$, and $\rho \gamma_\pm (\rho \otimes \rho) = \gamma_\mp$ and $(\rho \otimes \id) \fundsym = \fundsym (\id \otimes \rho)$). Alternatively, it is the symmetric monoidal subcategory of $\Cubez$ generated by $\Cubece$ and $\rho$, and it also follows from the work of \cite{mauri} that $\Cubecer$ is the free symmetric monoidal category on an involutive symmetric cubical monoid (where again $\fundsym$ is taken to be the symmetry isomorphism of the symmetric monoidal category). The category of \defterm{cubical sets with connections, symmetries, and reversals} is the presheaf category $\ccerset$.

There is also a subcategory $\Cubecr$ which has connections and reversals but not symmetries, and $\Cuber$ which has reversals but not connections or symmetries. These are considered in \cite{buchholtz-morehouse}. Note that if a cubical site has reversals and one connection, it automatically has both connections.

\textbf{Cubical sets with diagonals, connections, symmetries, and reversals.} Let $\Cubedcsr$ be the wide monoidal subcategory of $\Cubez$ generated by $\Cubecer$ along with the \defterm{diagonal maps} $\delta: \cube n \to \cube {2n}$, $x \mapsto (x,x)$. Then in fact we have $\Cubedcsr = \Cubez$. 
The diagonal is called a \defterm{contraction} by \cite{buchholtz-morehouse}.

There are also subcategories $\Cubeds$, $\Cubedsr$, $\Cubedcms$, $\Cubedcps$, $\Cubedcs$. Note that if a cubical site has diagonals, it automatically has symmetries, obtained by a diagonal followed by a projection: $(a,b) \mapsto (a,b,a,b) \mapsto (b,a)$.

\textbf{Conventions.} If $\Cubea$ is a cubical site and $X \in \cseta$, we write $X_n$ for $X_{\cube n} = X(\cube n)$. We write $\gamma_\pm$ to mean either $\gamma_\vee$ or $\gamma_\wedge$, with the understanding that in a given expression, the meaning of $\pm$ remains constant; we write $\gamma_\mp$ to mean $\gamma_\vee$ when $\gamma_\pm$ is $\gamma_\wedge$ and vice versa. We write $\iota^{i,\epsilon}$ for $\id \otimes \dots \otimes \iota^\epsilon \otimes \dots \otimes \id$ where the $\iota^\epsilon$ is in the $i$th slot, and similarly we write $\fundproj^i$, $\gamma_\pm^i$, $\rho^i$. We freely identify morphisms $s : \cube n \to \cube n$ built up from $\fundsym$ and its tensor composites with elements of the symmetric groups $\Sigma_n$. We sometimes also write $\iota^{i_1,\epsilon_1, \dots, i_k,\epsilon_k}$ to mean $\id \otimes \dots \otimes \iota^{\epsilon_1} \otimes \dots \otimes \iota^{\epsilon_k} \otimes \dots \otimes \id$, where $\iota^{\epsilon_j}$ is in the $i_j$th slot, and similarly we write $\fundproj^{i_1,\dots,i_k}$, $\rho^{i_1,\dots,i_k}$. We also write $\gamma_\pm$ for $\gamma_\pm(\gamma_\pm \otimes \id) \dots (\gamma_\pm \otimes \dots \otimes \id)$.


\textbf{Other cubical sites.} In \cite{buchholtz-morehouse}, cubical sites without projections (or \emph{weakenings} in their terminology) are also considered; presheaves on these are \emph{semi-cubical sets}. We will not consider such sites except in a few parenthetical remarks.


\section{Cubical sites without diagonals as Eilenberg-Zilber categories}\label{sec:cubsitegenreedy}

We work toward \cref{thm:cubereedy}, showing that all cubical sites without diagonals are strongly skeletal generalized Reedy categories. We are greatly aided by the fact that in the absence of diagonals, each output variable of a map $f : \cube n \to \cube m$ can depend on at most one input variable, as we make precise in \cref{def:depend}.

\begin{Def}
Let $\Cubea$ be a cubical site. Let $(\Cubea)_- \subset \Cubea$ be the monoidal subcategory generated by the projections and any connections, symmetries, and reversals present in $\Cubea$, and let $(\Cubea)_+$ be the monoidal subcategory generated by the face maps and any symmetries and reversals present in $\Cubea$.
\end{Def}

\begin{fact}\label{fact:fact}
Let $\Cubea$ be a cubical site without connections. By the results of \cite{grandis} (e.g. \cite[Equation (58)]{grandis}), we have that every morphism $f$ in $\Cubea$ may be factored in the form 
\[f = i c s r p\]
where 
\begin{itemize}
    \item $i$ is a composite of face maps;
    \item $c$ is a composite of connection maps in $\Cubea$;
    \item $s$ is a composite of symmetries in $\Cubea$;
    \item $r$ is a composite of reversals in $\Cubea$;
    \item $p$ is a composite of projections.
\end{itemize}
We may equally factor $f = i c r' s p$ where $r'$ is a (generally different) composite of reversals in $\Cubea$.
\end{fact}

\begin{Def}\label{def:depend}
Let $\Cubea$ be a cubical site, and let $f : \cube n \to \cube m$ be a morphism in $\Cubea$. For $j \leq n$ and $i \leq m$, we say that \defterm{the $i$th output variable of $f$ depends on its $j$th input variable} if there exists $(\epsilon_1, \dots, \epsilon_{n-1}) \in \{0,1\}^{n-1}$ and $j \leq n$ such that 
\[f_i(\epsilon_1, \dots, \epsilon_{j-1}, 0, \epsilon_j, \dots, \epsilon_{n-1}) \neq f_i(\epsilon_1, \dots, \epsilon_{j-1}, 1, \epsilon_j, \dots, \epsilon_{n-1}).\]
We say that $f$ \defterm{depends on all of its input variables} if, for every output variable $i$, there exists an input variable $j$ such that the $i$th output variable of $f$ depends on the $j$th input variable.

For $j$ an input variable, write for the moment $X_j = \{i \mid \text{the }i \text{th output variable depends on the }j\text{th input variable}\}$. If $j$ and If $j$ and $k$ are input variables such that $X_j = X_k$ is a subsingleton, then we say that $j,k$ are \defterm{codependent}.
\end{Def}

\begin{obs}\label{obs:depend-comp}
Let $\Cubea$ be a cubical site. 
\begin{enumerate}
    \item\label{obsitem:1} The $i$th output variable of $\id_{\cube n}$ depends on the $j$th input variable of $\id_{\cube n}$ if and only if $i = j$.
    \item\label{obsitem:2} Let $\cube p \xrightarrow g \cube n \xrightarrow f \cube m$ be a composite of maps therein. Let $k \in (0,p], i \in (0,m]$ be variables. Then the $i$th output variable of $fg$ depends on the $k$th input variable of $fg$ if and only if there exists $j \in (0,n]$ such that the $i$th output variable of $f$ depends in the $j$th input variable of $f$ and the $j$th output variable of $g$ depends on the $k$th input variable of $g$.
    \item\label{obsitem:3} Let $\cube {n'} \xrightarrow {f'} \cube{m'}$ be another map, and let $j \in (0,n+n'], i \in (0,m+m']$. Then the $i$th output variable of $f \otimes f'$ depends on the $j$th input variable of $f \otimes f'$ if and only if $i \leq m$ and $j \leq n$ and the ith output variable of $f$ depends on the $j$th input variable of $f'$, or $i > m$ and $j> n$ and the $i-m$th output variable of $f'$ depends on the $j-n$th input variable of $f'$.
    \item\label{obsitem:4} In particular, the class of maps depending on all of their input variables is closed under identities, composition and tensoring. It also contains the identities, so such maps span a wide monoidal subcategory of $\Cubea$.
\end{enumerate} 
\end{obs}


\begin{lem}\label{lem:depend}
Let $\Cubea$ be a cubical site, and let $c : \cube n \to \cube 1$ be a composite of connections, symmetries, and reversals in $\Cubea$. Then $c$ depends on all of its input variables. Moreover, for every input variable $j \leq n$, there exists a pseudo-section $\iota$ of $f$ in $\Cubea$ which is non-constant in the output variable $j$.
\end{lem}
\begin{proof}
Observe that the elementary connections, reversals, and symmetries depend on all of their input variables. So the first statement follows from \cref{obs:depend-comp}(\ref{obsitem:4}). For the second statement, if $f(\epsilon_1, \dots, -,\dots, \epsilon_{n-1})$ depends on the $j$th variable, then the face map $\iota$ which is constant at $\epsilon_{j'}$ for $j' \neq j$ and which is non-constant in output variable $j$, is a pseudo-section of $f$ (in that $f \iota$ is either the identity or a reversal; in the latter case, there must be some reversals appearing in $f$ (else $f$ would be order-preserving) and hence there must be reversals in $\Cubea$.)
\end{proof}

\begin{cor}\label{cor:codepend}
Let $\Cubea$ be a cubical site without diagonals, and let $f : \cube n \to \cube m$ be a degeneracy. Then codependency is an equivalence relation on the input variables of $f$. The associated partition is naturally ordered by the corresponding output variable. We call the associated (ordered) partition the \defterm{(ordered) dependency partition} associated to $f$, and the equivalence class of input variables upon which no output variable depends is the \defterm{trivial part}. Moreover, in the factorization $f = c r s  p$ of \cref{fact:fact}, the following hold:
\begin{enumerate}
    \item\label{item:2} The composite $c r$ decomposes as $c r = f_1 \otimes \cdots \otimes f_m$ for unique $f_i : \cube {n_i} \to \cube 1$ which are composites of connections and reversals in $\Cubea$. 
    The ordered dependency partition for $c r$ is of the form $\{\{1,2,\dots, n_1\},\{n_1+1,\dots, n_1+n_2\}, \dots, \{\dots,n_1 + \cdots + n_m\}\}$ with empty trivial part.
    \item\label{item:3} The ordered dependency partition for the composite $c r s$ is of the form $\{s\inv(\{1,2,\dots, n_1\}),s\inv(\{n_1+1,\dots, n_1+n_2\}), \dots, s\inv(\{\dots,m\})\}$ with empty trivial part.
    \item\label{item:4} The ordered dependency partition for $f = c r s  p$ differs from that for the composite $c r s$ by the interspersal of variables in the trivial part. The projection $p$ projects away precisely the variables in the trivial part.
\end{enumerate} 
\end{cor}
\begin{proof}
We first prove (\ref{item:2}). The decomposition $c r  = f_1 \otimes \cdots \otimes f_m$ and its uniqueness follows by definition, since the elementary connections and reversals all have codomain $\cube 1$ and everything else is built from these by tensoring and composition. 
The case $m = 1$ now follows from \cref{lem:depend}. By induction, the case of general $m$ then follows from \cref{obs:depend-comp}(\ref{obsitem:3}).
Thus (\ref{item:2}) holds, and in particular codependency is an equivalence relation when $p = s = \id_{\cube n}$.

From this, (\ref{item:3}) follows by \cref{obs:depend-comp}(\ref{obsitem:2}). In particular, we see that codependency is an equivalence relation when $p = \id_{\cube n}$.

For (\ref{item:4}), first observe that the case where $f = p$ is obvious (or can be pieced together from \cref{obs:depend-comp} by understanding the elementary projection), and the case where $p$ is trivial follows from (\ref{item:3}). The general case now follows by \cref{obs:depend-comp}(\ref{obsitem:2}).
\end{proof}

\begin{lem}\label{lem:depend+}
Let $\Cubea$ be a cubical site, and let $\cube n \xrightarrow f \cube m$ be a degeneracy therein. Then
\begin{enumerate}
    \item\label{lemitem:1} Suppose that $\iota : \cube m \to \cube n$ is a section of $f$ where $\iota$ is not constant in the output variables $j_1 < \dots < j_m$. Then for each $i \leq m$, we must have that the $i$th output variable of $f$ depends on the $j_i$th input variable of $f$.
    \item\label{lemitem:2} Suppose that $\iota : \cube m \to \cube n$ is a pseudo-section of $f$ where $\iota$ is not constant in the output variables $j_1 < \dots < j_m$.
    There exists a unique symmetry $s \in \Sigma_m$ such that for each $i \leq m$, the $i$th output variable of $f$ depends on the $j_{s\inv(i)}$th input variable of $f$. (If there are no symmetries in $\Cubea$, then $s$ is the identity.)
    \item\label{lemitem:3} Suppose that $1 \leq j < k \leq n$ are distinct input variables in the same part of the dependency partition. Then if $\iota$ is non-constant in the output variables $j$ and $k$, $\iota$ is not a pseudo-section of $f$.
    \item\label{lemitem:4} Suppose that $1 \leq j < k \leq n$ are distinct input variables in distinct nontrivial parts of the dependency partition. Then there exists a pseudo-section $\iota$ of $f$ which is non-constant in the output variables $j$ and $k$.
    \item\label{lemitem:5} Suppose that $1 \leq j \leq n$ is an input variable in the trivial part of the dependency partition. Then if $\iota$ is non-constant in the output variable $j$, then $\iota$ is not a pseudo-section of $f$.
    \item\label{lemitem:6} Suppose that $1 \leq j \leq n$ is an input variable in a nontrivial part of the dependency partition. Then there exists a pseudo-section $\iota$ of $f$ which is non-constant in the output variable $j$.
\end{enumerate}
\end{lem}
\begin{proof}
For (\ref{lemitem:1}), if the $i$th output variable of $f$ does not depend on the $j_i$th input variable of $f$, then $f\iota$ collapses the $i$th coordinate, so it is not an isomorphism.

For (\ref{lemitem:2}), if $\iota$ is a pseudo-section of $f$, then $\iota(f\iota)\inv$ is a section of $f$. Write $f\iota = r s$ where $r$ is composed of reversals in $\Cubea$ and $s$ is composed of symmetries in $\Cubea$. Then $s$ has the required properties by (\ref{lemitem:1}) and \cref{obs:depend-comp}(\ref{obsitem:2}).

For (\ref{lemitem:5}), observe that in this case $f\iota$ projects away the coordinate carried to the $j$th coordinate by $\iota_S$, and hence is not an isomorphism.

For (\ref{lemitem:6}), we decompose $f = crsp$ as in \cref{fact:fact}, and write $cr = f_1 \otimes \cdots \otimes f_m$ as in \cref{cor:codepend}(\ref{item:2}). Let $i$ be the output variable of $f$ upon which $j$ depends. Then $p(s(j)) - n_1 - \cdots - n_{i-1}$ is a in input variable upon which $f_i$ depends. By \cref{lem:depend}, there exists a pseudo-section $\iota_i$ of $f_i$ depending on this input variable. Choose arbitrary pseudo-sections $\iota_{i'}$ for the other $f_{i'}$'s. Then $\iota_1 \otimes \otimes \iota_m$ is a pseudo-section of $cr$, and $s\inv(\iota_1 \otimes \cdots \otimes \iota_m)$ is a pseudo-section for $crs$. Composing with an arbitrary section of $p$, we obtain a pseudo-section of $f$ with the desired properties.

For (\ref{lemitem:3}), observe that in this case $f \iota$ must collapse variable $j,k$ down to one variable, so it is not an isomorphism.

For (\ref{lemitem:4}), we reason as for (\ref{lemitem:6}) except that now both $\iota_i$ and $\iota_h$ are chosen non-arbitrarily (where $i$ depends on $j$ and $h$ depends on $k$).
\end{proof}

\begin{cor}\label{cor:codepend+}
Let $\Cubea$ be a cubical site, and let $f$ be a degeneracy therein. Then the unordered dependency partition of $f$ and the trivial part thereof are determined by the set of pseudo-sections of $f$. If there are no symmetries in $\Cubea$, then the ordered dependency partition of $f$ is determined by the set of pseudo-sections of $f$.
\end{cor}
\begin{proof}
By \cref{lem:depend+}(\ref{lemitem:3} and \ref{lemitem:4}), we know exactly which input variables are in the trivial part, namely the variables in which every pseudo-section is constant. By \cref{lem:depend+}(\ref{lemitem:5} and \ref{lemitem:6}), we then know exactly which pairs of nontrivial input variables are codependent. Thus we recover the unordered dependency partition and its trivial part. If there are no symmetries in $\Cubea$, then we may deduce which part corresponds to which output variable because we know that the nontrivial parts are intervals ordered as in the output variables by \cref{cor:codepend}(\ref{item:3}).
\end{proof}

We thank Denis-Charles Cisinski for suggesting the proof idea of the following theorem:
\begin{thm}[cf. \cite{isaacson}]\label{thm:cubereedy}
Let $\Cubea$ be a cubical site without diagonals. Then $(\Cubea, (\Cubea)_-,(\Cubea)_+)$ is a strongly skeletal generalized Reedy category, and in particular an Eilenberg-Zilber category. If $\Cubea$ has no symmetries or reversals, then it is even a regular skeletal Reedy category.
\end{thm}
\begin{proof}
The fact that $\Cubea$ is a generalized Reedy category (where $\deg(\cube n) = n$) follows from the $icsrp$ factorization of \cref{fact:fact}. As $i \in (\Cubea)_+$ and $csrp \in (\Cubea)_-$, this gives existence of factorizations. For uniqueness of factorizations, first observe that if $\phi$ is an automorphism in $\Cubea$ and $i$ is a composite of face maps, then $\phi i = i \psi$ for another automorphism $\psi$ in $\Cubea$. So if $f = iq = i'q'$ with $i,i' \in (\Cubea)_+$, then we may assume that $i,i'$ are composites of face maps. Then we must have $i = i'$ because $i$ is determined by the image of $f$. Since composites of face maps are monic, it follows that $q = q'$ as well.

Now let us show that $\Cubea$ is strongly skeletal. We begin by observing that any morphism of $(\Cubea)_-$ is a split epimorphism. We must show that every split epimorphism $\psi : \cube n \to \cube m$ of $(\Cubea)_-$ is determined up to pseudo-equality by its set of pseudo-sections. We argue in two steps:

\begin{component}[Case 1: $m = 1$]
The pseudo-sections of $\psi$ pseudo-uniquely determine $\psi$ when $m=1$. First consider the case when $\Cubea$ has no reversals. In this case, pseudo-sections are the same as sections, and pseudo-equality is the same as equality. Observe that $\psi: \cube n \to \cube 1$ is a Boolean function, uniquely determined by the set $A$ of vertices mapped to $0$ and the complementary set $B$ of vertices mapped to $1$ (these sets are nonempty because $\psi$ is surjective). Let $E$ be the set of edges in $\cube n$ with one vertex in $A$ and the other in $B$. Because $\psi$ is surjective, both $A$ and $B$ are nonempty. If $v \in A$, then the neighboring vertices are in $A$ or $B$ according as the edge between them is not or is in $E$. Since the 1-skeleton of $\cube n$ is connected, we may start with $v_0 \in A$, and then progress outward through the graph, reading off exactly which vertices are in $A$ and which are in $B$, so we determine $\psi$ from its set of sections.

Now consider the case where $\Cubea$ has reversals. Define $A,B$ as before. Then $\psi$ is determined by the partition $\{A,B\}$ up to pseudo-equality. By the same procedure as before, we may start with some vertex $v_0$ and determine exactly which vertices are in the same part of the partition as $v_0$ and which are in the other part by looking at $E$. Thus we have determined $A$ and $B$ up to swapping of factors, and thus up to possibly passing from $\psi$ to the pseudo-equal $\rho \psi$, as desired.
\end{component}

\begin{component}[Case 2: General Case]
By \cref{cor:codepend+}, the set of pseudo-sections of $\psi$ determines the dependency partition of $\psi$ and its trivial part. By \cref{cor:codepend}, the trivial part determines $p$. Composing with a symmetry in $\Cubea$, we may pass to a pseudo-equal $\psi$ where the dependency partition consists of intervals which do not intersperse (except for the trivial part). In this case, by \cref{cor:codepend}(\ref{item:3}), we have that in the factorization $\psi = crsp$ of \cref{fact:fact}, we have $s = s_1 \otimes \cdots \otimes s_m$, where $s_j : \cube {n_j} \to \cube {n_j}$ is a symmetry in $\Cubea$. Thus we are reduced to the case where there exists a tensor splitting $\psi = \psi_1 \otimes \cdots \otimes \psi_m$ where $\psi_j : \cube{n_j} \to \cube 1$. In this case, $\iota$ is a pseudo-section of $\psi_1 \otimes \cdots \otimes \psi_m$ if and only if $\iota = (\iota_1 \otimes \cdots \otimes \iota_m)$, where each $\iota_j$ is a pseudo-section of $\psi_j$, and $a$ is an automorphism in $\Cubea$. So if $\psi = \psi_1 \otimes \cdots \otimes \psi_m$ has the same set of pseduo-sections as $\psi' = \psi'_1 \otimes \cdots \otimes \psi'_m$, it follows that there is a permutation $s'$ lying in $\Cubea$ such that $\psi_j$ has the same set of pseudo-sections as $\psi'_j$. By the first part, $\psi_j$ is pseudo-equal to $\psi'_{s'(j)}$, i.e. $\psi'_{s'(j)} = b_j\psi_j$ where $b_j$ is an automorphism in $\Cubea$. Thus, setting $b = b_1 \otimes \cdots \otimes b_m$, we have that $\psi' = ba\psi$, so that $\psi$ and $\psi'$ are pseudo-equal as desired.
\end{component}

Thus $\Cubea$ is a strongly skeletal generalized Reedy category. If $\Cubea$ has no symmetries or reversals, then of course $\Cubea$ is a skeletal Reedy category. It is regular because the face maps are monomorphisms.
\end{proof}

This allows us to recover and generalize a result of Isaacson, without directly verifying absolute pushout diagrams:
\begin{cor}[\cite{isaacson} for $\Cubea = \Cubecme, \Cubecpe$]
Let $\Cubea$ be a cubical site without diagonals. Then $(\Cubea,(\Cubea)_-,(\Cubea)_+)$ is an Eilenberg-Zilber category.
\end{cor}
\begin{proof}
This follows from \cref{thm:cubereedy} and \cref{thm:ssgrEZ}.
\end{proof}

\section{Cubical sites with diagonals}\label{sec:diag}

To wrap things up, we take a look at cubical sites with diagonals. Some are EZ and others are not (\cref{thm:diag}). In order to carry out our analysis, we go through several descriptions of these sites, some of which are folklore.

\begin{prop}
$\Cubedcms$ (resp. $\Cubedcps$) is canonically identified with the category of finite Boolean algebras and binary-meet-preserving (resp. binary-join-preserving) maps.
\end{prop}
Note that the maps of $\Cubedcms$ do \emph{not} typically preserve the top element $\top$. Dually, the maps of $\Cubedcps$ do \emph{not} typically preserve the bottom element $\bot$.
\begin{proof}
We will prove the case of $\Cubedcms$; the case of $\Cubedcps$ is dual. We identify $\cube n$ with the Boolean algebra $2^n$ given as the powerset of $\{1,\dots,n\}$ in the natural way, where a vertex $(v_1,\dots,v_n)$ corresponds to the subset $\{i \mid v_i = 1\} \subseteq \{1,\dots,n\}$. Under this identification, the morphisms of $\Cubedcms$ are binary-meet-preserving maps; to see this it suffices to check on generators. This is obvious for projections and diagonals. For $\gamma_-$ and $\iota^1$, it follows because limits commute with limits, and for $\iota^0$ it holds because constants commute with nonempty limits.

Conversely, given a binary-meet-preserving map $f: 2^m \to 2^n$, we wish to factor $f$ into maps built from the generators of $\Cubedcms$. Because the inclusion of $\Cubedcms$ into binary-meet-semilattices is strong cartesian monoidal, we may reduce to the case $n=1$. For $m=0$, the result is clear. By binary-meet-preservation, for $m \geq 1$ we have $f(v_1,\dots,v_m) = f((v_1,1,\dots,1) \wedge \cdots \wedge (1,\dots,1,v_m)) = f(v_1,\dots,1)\wedge \cdots \wedge f(1,\dots,1,v_m)$, which allows us to reduce to the case $m=1$, which is also clear.
\end{proof}

\begin{prop}
$\Cubedcms$ (resp. $\Cubedcps$) is the free cartesian monoidal category on a half-cubical monoid.
\end{prop}
\begin{proof}
Again we treat $\Cubedcms$, as $\Cubedcps$ is dual. In \cite{buchholtz-morehouse}, it is shown that the cartesian monoidal theory of $\{0,1\}$ in this signature coincides with the theory of $[0,1]$ in this signature. An axiomatization is given. 
By definition, the theory of $\{0,1\}$ in this signature comprises precisely the equations satisfied in our category $\Cubedcms$. The axioms given are precisely the axioms of a half-cubical monoid.
\end{proof}

We would like to thank Keith Kearnes for indicating the proof of the following:
\begin{prop}\label{prop:msl}
$\Cubedcms$ (resp. $\Cubedcps$) does not have split idempotents. Its Karoubi envelope is canonically identified with the category of finite distributive lattices and binary-meet semilattice maps (resp. binary-join semilattice maps).
\end{prop}
\begin{proof}
Again we treat $\Cubedcms$, as $\Cubedcps$ is dual. Recall that Birkoff's Representation Theorem \cite{birkhoff} implies that a finite lattice is distributive if and only if it is isomorphic to the lattice of down-closed sets in a finite poset $P$. For any such lattice $2^P$, the forgetful functor $2^P \to 2^{|P|}$ has both left and right adjoints, where $|P|$ is the underlying set of the poset $P$. By composing with the right adjoint, we exhibit $2^P$ as a retract of the Boolean algebra $2^{|P|}$ via meet-preserving maps.

Conversely, \cite{horn-kimura} characterizes the finite distributive lattices as the finite injective binary-meet-semilattices, so finite distributive lattices are closed under retracts.
\end{proof}

\begin{prop}
$\Cubeds$ is the free cartesian monoidal category on a bipointed object.
\end{prop}
\begin{proof}
In \cite{buchholtz-morehouse}, it is shown that the cartesian monoidal theory of $\{0,1\}$ in this signature coincides with the theory of $[0,1]$ in this signature. An axiomatization is given. 
By definition, the theory of $\{0,1\}$ in this signature comprises precisely the equations satisfied in our category $\Cubeds$. The axioms given are precisely the axioms of an object equipped with two points $0,1$.
\end{proof}

\begin{prop}
$\Cubedsr$ is the free cartesian monoidal category on a pointed object with an involution.
\end{prop}
\begin{proof}
In \cite{buchholtz-morehouse}, it is shown that the cartesian monoidal theory of $\{0,1\}$ in this signature coincides with the theory of $[0,1]$ in this signature. An axiomatization is given.
By definition, the theory of $\{0,1\}$ in this signature comprises precisely the equations satisfied in our category $\Cubedsr$. The axioms given are precisely the axioms of an object equipped with a point $0$ and an involution.
\end{proof}

\begin{prop}
$\Cubedcs$ is the free cartesian monoidal category on a cubical monoid.
\end{prop}
\begin{proof}
In \cite{buchholtz-morehouse}, it is shown that the cartesian monoidal theory of $\{0,1\}$ in this signature coincides with the theory of $[0,1]$ in this signature. An axiomatization is given.
By definition, the theory of $\{0,1\}$ in this signature comprises precisely the equations satisfied in our category $\Cubedcsr$. The axioms given are precisely the axioms of an object equipped with a cubical monoid structure.
\end{proof}

\begin{prop}[Kapulkin]\label{prop:kap}
$\Cubedcs$ is canonically identified with the category of finite Boolean algebras and order-preserving maps.
\end{prop}
\begin{proof}
Clearly every map of $\Cubedcs$ is order-preserving. Conversely, if $f: \cube n \to \cube 1$ is order-preserving, then we may write $f(x_1,\dots, x_n) = f(x_1,\dots,x_{n-1},0) \vee (f(x_1,\dots,x_{n-1},1) \wedge x_n)$. The $(n-1)$-ary function $f(x_1,\dots,x_{n-1},\varepsilon)$ is order-preserving for $\varepsilon = 0,1$, and so lies in $\Cubedcs$ by induction on $n$. Thus $f$ also lies in $\Cubedcs$.
\end{proof}

\begin{prop}[Sattler]\label{prop:satt}
$\Cubedcs$ does not have split idempotents. Its Karoubi envelope is canonically identified with the category of finite lattices and order-preserving maps.
\end{prop}
\begin{proof}
By the Knaster-Tarski fixed point theorem \cite{knaster}, every retract of an order-preserving idempotent on a lattice is finite. So the idempotent completion of $\Cubedcs$ is contained in the category of finite lattices. Conversely, every finite lattice embeds via Yoneda into its lattice of down-closed sets, and the right adjoint is an order-preserving retract. The lattice of down-closed sets in turn embeds into the boolean algebra of all subsets, and the left adjoint is an order-preserving retract. Thus every finite lattice is an order-preserving retract of a finite boolean algebra as desired.
\end{proof}

\begin{prop}\label{prop:bool}
$\Cubedcsr$ is the free cartesian monoidal category on a Boolean algebra.
\end{prop}
\begin{proof}
Clearly $\cube 1$ is a Boolean algebra. Moreover, $\cube n$ is the free Boolean algebra on $n$ generators. Thus if two morphisms $\cube n \rightrightarrows \cube 1$ constructed from Boolean operations are equal, then this is implied by the theory of a Boolean algebra. So there are no further relations which hold in $\Cubedcsr$.
\end{proof}

\begin{prop}\label{prop:cubez}
$\Cubedcsr$ is equivalent to the category $\Cubez$ of finite powersets and all set-maps.
\end{prop}
\begin{proof}
By using arbitrary Boolean operations, it is clear that we can construct all maps of $\Cubez$.
\end{proof}

\begin{prop}\label{prop:fin}
$\Cubedcsr$ does not have split idempotents. Its Karoubi envelope is equivalent to the category of finite nonempty sets.
\end{prop}
\begin{proof}
This follows from Proposition \ref{prop:cubez}
\end{proof}

We now give the main theorem of this section.

\begin{thm}\label{thm:diag}
\begin{enumerate}
    \item The cubical sites $\{\Cubeds,\Cubedsr\}$ are strongly skeletal generalized Reedy categories with the obvious degree functions, and in particular Eilenberg-Zilber categories.
    \item The cubical sites $\{\Cubedcms,\Cubedcps,\Cubedcs\}$ are not idempotent complete (and hence not generalized Reedy categories). Even after idempotent completion, the split epimorphisms do not form the left half of an orthogonal factorization system on the idempotent-completed sites, and in particular the idempotent-completed sites are not Eilenberg-Zilber categories for any degree function.
    \item The cubical site $\Cubedcsr$ is not idempotent complete (and hence not a generalized Reedy category). But after idempotent completion, the idempotent-complete site is the category of finite nonempty sets, and is a strongly skeletal generalized Reedy with the obvious degree function, and in particular Eilenberg-Zilber.
\end{enumerate} 
\end{thm}
\begin{proof}
$\Cubeds$ is the free cartesian monoidal category on a bipointed object (Proposition \cite{cisinski} 8.4.5, \cite{grandis}). The category of free $\Set$-models of finite rank for this Lawvere theory is the category $\Fin_{\ast\ast}$ of finite sets equipped with an ordered pair of distinct basepoints $(\ast_0,\ast_1)$. Thus $\Cubeds \simeq \Fin_{\ast\ast}^\op$ is equivalent to the opposite category. In $\Fin_{\ast\ast}$, there is an orthogonal (Mono,Epi) factorization system, and every monomorphism splits. Two monomorphisms are pseudo-equal iff they have the same image, and a pseudo-retract is precisely a map which restricts to a bijection on this image. Thus two monomorphisms are pseudo-equal iff they have the same set of pseudo-retracts, and $\Cubeds$ is strongly skeletal.

$\Cubedsr$ is the free cartesian monoidal category on an object equipped with an involution and a point (\cite{grandis}). The category of free $\Set$-models of finite rank for this Lawvere theory is the category $\Fin^{C_2,\textrm{free}}_{C_2/}$ of finite free $C_2$-sets equipped with an equivariant map from $C_2$. Thus $\Cubedsr \simeq (\Fin^{C_2,\textrm{free}}_{C_2/})^\op$ is the opposite category. In $\Fin^{C_2,\textrm{free}}_{C_2/}$, there is an orthogonal (Mono,Epi) factorization system (since free $C_2$-sets and finite $C_2$-sets are closed under subobjects), and every monomorphism splits. The description of pseudo-retracts and pseudo-equality is similar to that in $\Fin_{\ast\ast}$, and thus $\Cubedsr$ is strongly skeletal.

The cubical sites $\{\Cubedcms,\Cubedcps\}$ are equivalent; we treat the latter. It is equivalent to the category of finite Boolean algebras with binary-meet-semilattice homomorphisms as morphisms. This category is not idempotent complete; its idempotent splitting is the category $\FinDist^\wedge$ of finite distributive lattices with binary-meet-semilattice homomorphisms (Proposition \ref{prop:msl}). Note that all surjections of $\FinDist^\wedge$ split: For finite Boolean algebras are free semilattices, and it is easy to see that this implies they are projective among binary-meet-semilattices, and every object of $\FinDist^\wedge$ is a retract of one of these. Any morphism $X \to Y$ in $\FinDist^\wedge$ right orthogonal to the split epimorphisms is in particular right orthogonal to the projection $\cube 1 \to \cube 0$; because lattices are connected this implies that $X \to Y$ is monic. Conversely, the split epimorphisms are left orthogonal to every monic morphism. So if $\FinDist^\wedge$ were Eilenberg-Zilber, it would have a split-epi / mono factorization system. But let $P$ be a finite meet-semilattice which is not distributive. Let $D \to P$ be a surjective meet-semilattice homomorphism and $P \to 2^P$ the Yoneda embedding. Then the composite map $D \to P \to 2^P$ is a morphism in $\FinDist^\wedge$ which we claim does not admit a split-epi / mono factorization. For its surjection / mono factorization in the category of meet-semilattices is given by $D \to P \to 2^P$, and does not factor through a distributive lattice; a split-epi / mono factorization in $\FinDist^\wedge$ would necessarily be a \emph{different} surjection / mono factorization in the category of meet-semilattices, violating the uniqueness of such factorizations.

The cubical site $\Cubedcs$ is equivalent to the category of finite Boolean algebras with order-preserving maps as morphisms (Proposition \ref{prop:kap}). This category is not idempotent complete; its idempotent splitting is the category $\FinLat^\leq$ of finite lattices with order-preserving maps (Proposition \ref{prop:satt}). As before, the right orthogonal complement of the split epimorphisms comprises the monomorphisms. Let $W$ be the poset $y_1 > p < x > q < y_2$, and let $P = 2^W / (y_1 \sim  y_2)$ where we identify elements of $W$ with their images under the downset (Yoneda) embedding $W \to 2^W$, and again with their images in $P$. Let also $y$ denote the common image of $y_1,y_2$ in $P$. Then $P$ is not a lattice. For if it were, then because $p,q \leq x,y$ in $P$, there would be $a \in P$ such that $p,q \leq a \leq x,y$. We can't have $a = y$ since $y \not \leq x$, so $a$ has a unique lift to $2^W$. Therefore $p,q \leq a$ implies $p \vee q \leq a$ in $2^W$. But $p \vee q \not \leq y_1$ and $p \vee q \not \leq y_2$ in $2^W$; it follows that $p \vee q \not \leq y$ in $P$. If there were such an $a$, we would have $p \vee q \leq a \leq y$, a contradiction. Thus $P$ is not a lattice. Now, consider the composite map $2^W \to P \to 2^P$ with the Yoneda embedding. This is an order-preserving map of finite lattices whose quotient / injection factorization is given as $2^W \to P \to 2^P$; if there were also a split-epi / injection factorization through a lattice, this would violate the uniqueness of quotient / injection factorizations of posets. Hence no split epi / injection factorization of this map exists, and so $\FinLat^\leq$ is not Eilenberg-Zilber.

The cubical site $\Cubedcsr$ is equivalent to the category of finite Boolean algebras with all set maps as morphisms (Proposition \ref{prop:cubez}). This category is not idempotent complete; its idempotent splitting is the category $\Fin_{\neq \emptyset}$ of finite nonempty sets with set maps between them (Proposition \ref{prop:fin}). This category has unique (split epi, mono) factorizations, and is EZ with respect to the obvious degree function.
\end{proof}

\bibliographystyle{alpha}
\bibliography{ez-with-cubes}

\end{document}